\documentclass[12pt]{amsart}
\usepackage{latexsym, ifthen,amsfonts,amsmath,amssymb}
\usepackage[dvips]{graphics}
\usepackage{epsfig,color}
\usepackage[margin=1.5in]{geometry}

 \newtheorem{theorem}{Theorem}

 \newtheorem{lemma}[theorem]{Lemma}
 \newtheorem{corr}[theorem]{Corollary}
 
 \newtheorem{proposition}[theorem]{Proposition}
 \newtheorem{deff}[theorem]{Definition}
 
 \newcommand{\bth}{\begin{theorem}}
 \newcommand{\ble}{\begin{lemma}}
 \newcommand{\bcor}{\begin{corr}}
 \newcommand{\bdeff}{\begin{deff}}
 \newcommand{\bprop}{\begin{proposition}}
 \renewcommand{\eth}{\end{theorem}}
 \newcommand{\ele}{\end{lemma}}
 \newcommand{\ecor}{\end{corr}}
 \newcommand{\edeff}{\end{deff}}
 \newcommand{\beq}{\begin{equation}}
 \newcommand{\eeq}{\end{equation}}
 \newcommand{\eprop}{\end{proposition}}

 \renewcommand{\Pi}{\varPi}

 \renewcommand{\epsilon}{\varepsilon}

\renewcommand{\S}{\Sigma}

\begin{document}
\begin{title}
{Curvature estimates for minimal surfaces with total boundary curvature less than 4$\pi$.}
\end{title}
\begin{author}
{Giuseppe Tinaglia}
\end{author}
\date{\today}
\maketitle

\date{\today}
\begin{abstract}We establish a curvature estimate for classical minimal surfaces with total boundary curvature less than 4$\pi$. The main application is a bound on the genus of these surfaces depending solely on the geometry of the boundary curve. We also prove that the set of simple closed curves with total curvature less than $4\pi$ and which do not bound an orientable compact embedded minimal surface of genus greater than $g$, for any given $g$, is open in the $C^{2,\alpha}$ topology.
\end{abstract}

\section*{Introduction}
\pagestyle{plain}

In this paper we discuss the geometry and topology of compact minimal surfaces whose boundary $\Gamma$ has total curvature $T(\Gamma)$ less than $4\pi$; Here $T(\Gamma)=\int_\Gamma|k|$,  where $k$ is the curvature of $\Gamma$. If $\Gamma$ is a connected simple closed curve, we denote by $M(\Gamma)$ the family of all orientable compact minimal surfaces whose boundary is $\Gamma$. It is well known that if $\Gamma$ is a connected piecewise $C^1$ simple closed curve, then $M(\Gamma)$ contains at least an immersed minimal disk \cite{do3,os7,ra1}. We stress that Ekholm, White, and Wienholtz proved in \cite{eww1} that a classical minimal surface with boundary $\Gamma$ such that $T(\Gamma)\leq 4\pi$ must in fact be embedded, regardless of the topological type.

The main theorem of this paper is the following curvature estimate:
\begin{theorem}\label{mainthm}
Let $\Gamma\subset \mathbb{R}^3$ be a $C^{2,\alpha}$ connected simple closed curve such that the total curvature 
$T(\Gamma)$ is strictly less than $4\pi$. Then there exists a constant $C=C(\Gamma)$ such that
$$\sup_{\Sigma \in{M(\Gamma)}}|K_\Sigma|\leq C,$$
where $K_\Sigma$ is the Gaussian curvature of $\Sigma$.

The constant $C$ depends on the radius of the largest embedded tubular neighborhood around $\Gamma$, as well as on upper bounds for $\|\Gamma\|_{C^{2,\alpha}}$ and the length of $\Gamma$.\end{theorem}

A finer quantative version of Theorem~\ref{mainthm} is to be found in Section~\ref{constant}.

The proof of Theorem~\ref{mainthm} is based on a compactness argument, which we now sketch. Suppose that there exists a sequence of embedded minimal surfaces and points on these surfaces which remain away from the boundary, where the curvature blows up. This being the case, we use a rescaling argument to obtain a new sequence of surfaces that converges to a complete nonplanar embedded minimal surface. A key point in the argument is a delicate application of the density estimates in \cite{eww1} which are used to prove that the limiting surface must in fact be a plane, thus yielding a contradiction. A similar, yet more refined argument, is needed when the sequence of points converges to a boundary point.

Using the Gauss-Bonnet formula, together with an area estimate for minimal surfaces, our curvature estimate implies a bound on the genus of each $\Sigma$ in $M(\Gamma)$ that depends only on the geometry of $\Gamma$. Consequently, the topology of the elements in $M(\Gamma)$ cannot be arbitrary. In addition to Theorem~\ref{mainthm}, we prove that the bound on the Gaussian curvature varies continuously with $\Gamma$, relative to the $C^{2,\alpha}$ topology. This is used to establish that {\it the set of all simple closed curves with $T(\Gamma)<4\pi$, and which \underline{do} \underline{not} bound an embedded and orientable minimal surface of genus greater than $g$, for any given $g$, is open in the $C^{2,\alpha}$ topology} (Theorem~\ref{maincor}). 

There are many special curves $\Gamma$ with total curvature less than $4\pi$ for which it has been proved that $M(\Gamma)$ consists only of disks (e.g., if $\Gamma$ lies on the boundary of a convex set \cite{my4b}). Theorem~\ref{maincor} reveals that if $\Gamma'$ is a new curve which is obtained by a slight modification of such a $\Gamma$, then $M(\Gamma')$ must still consist only of disks. These results relate to the Ekholm-White-Wienholtz conjecture on the nature of minimal surfaces spanning a curve $\Gamma$ with $T(\Gamma)\leq 4\pi$ (see \cite{eww1}).

This paper is organized as follow. In the first section we prove Theorem~\ref{mainthm}. Next we examine the constant in Theorem~\ref{mainthm} in order to prove that it depends only on a few geometric quantities associated with the curve $\Gamma$. Finally, we discuss some interesting applications of our curvature estimate.

\section*{Acknowledgments}
It is a pleasure to thank Professor Brian White for the many insightful conversations we had while the author was a post-doctoral fellow at Stanford University.

\section{Proof of Theorem \ref{mainthm}}\label{connected}

We begin by stating the density estimate in~\cite{eww1} which is needed in our arguments. In the statement below $\Theta(M,p)$ stands for the density of $M$ at $p$ --- see \cite{eww1} for details.
\begin{proposition}
Given $\theta>0$ there exists $\delta=\delta(\theta)>0$ such that the following holds. Let $\Gamma$ be a simple closed curve with $T(\Gamma)\leq \theta< 4\pi$ and let $p\in \Sigma \in M(\Gamma)$. Then,
\begin{equation}\label{est1}
\Theta(M,p)\leq 2-\delta, \quad \quad \text{if} \quad p\in M(\Gamma)\backslash \Gamma
\end{equation}
\begin{equation}\label{est2}
\Theta(M,p)\leq \frac{3}{2} -\delta, \quad \quad \text{if} \quad p\in \Gamma.
\end{equation}
\end{proposition}

Notice that the fact that $\Gamma$ is $C^{2,\alpha}$ implies that for any $\S\in M(\Gamma)$ there exists a constant $C(\S)$ which bounds the Gaussian curvature (see~\cite{ni3b}). In this paper, we are showing that the bound does not depend on the surface, but only on the geometry of the boundary. We also note that for a minimal surface, $-2K_{\Sigma}=|A|^2$. Here, as usual, if we denote by $k_1$ and $k_2$ the principal curvatures, $|A|=\sqrt{k_1^2+k_2^2}$ is the norm of the second fundamental form, and $K_{\Sigma}=k_1k_2$ is the Gaussian curvature.\\

\noindent{\it Proof of Theorem \ref{mainthm}.} To prove the bound on the Gaussian curvature we use a compactness argument. Assuming that the statement is false, we can find a sequence of minimal surfaces $\S_n$ such that $$\max_{\S_n}|A_n|> n.$$ Let $p_n \in \S_n$ such that $$|A_n(p_n)|=\max_{\S_n}|A_n|> n,$$ and let $q_n\in \Gamma$ such that 
$$r_n=\text{dist}_{\S_n}(\Gamma, p_n)= \text{dist}_{\S_n}(q_n,p_n).$$\\

\noindent\textsl{Alternative I: $r_n|A_n(p_n)|$ tends to infinity.} Consider the sequence of minimal surfaces, which we will still call $\S_n$, obtained by rescaling the connected component of $\S_n\cap B_{r_n}(p_n)$ which contains $p_n$, by a factor $|A(p_n)|$; Here $B_r(p)$ denotes the ball of radius $r$ centered at $p$. After a translation that takes $p_n$ to the origin, $\S_n$ becomes a sequence of embedded minimal surfaces with $|A_n|$ uniformly bounded by one, and $|A_n(0)|=1$. After possibly going to a subsequence, since $\lim_{n\rightarrow \infty} r_n|A_n(p_n)|=\infty$, we can assume that $\S_n$ converges to an orientable properly embedded minimal surface $\overline{\S}$ whose norm of the second fundamental form is bounded and has value one at the origin (in particular, $\overline{\S}$ is not a plane). Thanks to the results in \cite{eww1} there exists $\delta>0$ such that

\begin{equation}\label{density1}
\frac{Area(\S_n\cap B_r(q))}{\pi r^2}<2-\delta<2, \quad \text{ for any } q\in \S_n \text{ and any } r>0.
\end{equation}

In particular, \eqref{density1} gives

\begin{equation}\label{density2}
\frac{Area(\overline{\S}\cap B_r(q))}{\pi r^2}<2-\delta<2, \quad \text{ for any } q\in \overline{\S} \text{ and any } r>0.
\end{equation}

Our goal is to reach a contradiction by showing that $\overline{\S}$ must be a plane. Since $\overline{\S}$ has quadratic area growth, it is possible to take its cone at infinity. In other words, there exists a sequence $t_n>0$ approaching zero such that $t_n\overline{\S}$ converges to a stationary cone $C$ whose density at the origin must be less than $2-\delta$. It is known that the intersection of $C$ with the unit sphere consists of a collection of geodesic arcs. We claim that the said intersection consists of a single great circle with multiplicity one, which would imply that $\overline{\S}$ is a plane. In order to prove the claim we rule out the possibility that there is a point which is the end point of more than two arcs. Clearly, there cannot be a point where more than three arcs meet, otherwise the density at that point would be at least 2, and that would contradict \eqref{density2}. If there exists a point where exactly three arcs meet, then we can find a large circle which is transversal to $\overline{\S}$ and intersects $\overline{\S}$ in exactly three points. However, this contradicts the fact that the intersection number of the circle and the minimal surface is zero mod 2 (see~\cite{gupo}).\\

\noindent\textsl{Alternative II: $r_n|A_n(p_n)|$ is bounded.} Note that $|p_n-q_n||A_n(p_n)|\leq r_n|A_n(p_n)|$ is bounded. Reasoning as in the previous case, we obtain a sequence of minimal surfaces $\S_n$ that converges to an orientable properly embedded minimal surface $\overline{\S}$, bounded by a straight line $L$. As before, $\overline{\S}$ has bounded second fundamental form with lenght one at a certain point (in particular, $\overline{\S}$ is not a half-plane). Moreover, there exists $\delta>0$ such that

\begin{equation}\label{density1a}
\frac{Area(\overline{\S}\cap B_r(q))}{\pi r^2}<2-\delta<2, \quad \text{ for any } q\in M\backslash L \text{ and any } r>0,
\end{equation}
and 
\begin{equation}\label{density2a}
\frac{Area(\overline{\S}\cap B_r(q))}{\pi r^2}<\frac{3}{2}-\delta<\frac{3}{2}, \quad \text{ for any } q\in L \text{ and any } r>0.
\end{equation}

Our goal is to reach a contradiction by showing that $\overline{\S}$ must be a half-plane. Let $C$ be the cone at infinity. The intersection of $C$ with the unit sphere consists of a collection of geodesic arcs which must contain two antipodal points $P, Q$. We can always assume $P=(0,1,0)$, $Q=(0,-1,0)$. We claim that the intersection of $C$ with the unit sphere consists of a single half great circle. This implies that $\overline{\S}$ is a half-plane.  Similarly to the previous case, one can show that the intersection of $C$ with the unit sphere consists of halves of great circles having $P, Q$ as the endpoints. Moreover, there can be at most two great circles otherwise the density at the origin would be at least $\frac{3}{2}$, and that contradicts \eqref{density2a}. Assume therefore that there are two half great circles. Using the Schwarz reflection principle one obtains an orientable properly immersed minimal surface $M'$ without boundary. Let $C'$ be its cone at infinity. We have already shown that the intersection of $C'$ with the unit sphere consists of at most two great circles. Let $\pi$ be a plane through the origin which does not contain either great circle, and let 
\begin{multline}\omega(R)=\{(0,0,t): -R\leq t\leq R\}\cup \{(t,0,R): 0\leq t\leq R\}\cup\\
\cup \{(R,0,t): -R\leq t\leq R\}\cup \{(t,0,-R): 0\leq t\leq R\}.\end{multline} If there are two great circles, because of the simmetry of $M'$, for large $R$ the union of curves, $\omega(R)$, intersects $M'$ in an odd number of points. But, again, this contradicts the fact that the mod 2 intersection number of $\omega(R)$ and the minimal surface is zero.
\qed\\

\noindent{\bf Remarks}
\begin{enumerate}
\item The conlusion in Theorem~\ref{mainthm} may fail if $\Gamma$ is only assumed to be continuous. Certainly there exist continuous (but not $C^{2,\alpha}$) simple closed curves $\Gamma$'s and minimal surfaces $\S\in M(\Gamma)$ such that $\sup_{\S}|K_\S|=\infty$ (see \cite{ni3b,ni2}).

\item The connectedness hypothesis in Theorem~\ref{mainthm} is redundant and was included only to preserve the flow of the presentation. In fact, by the Fenchel-Borsuk theorem the total curvature of a connected simple closed curve $\Gamma$ is always greater than or equal to $2\pi$, with equality holding if and only if $\Gamma$ is a convex planar curve (see \cite{bors1,fa1,fe1,mil1}).
\end{enumerate}

\section{The constant C}\label{constant}

A slight modification of the argument in the proof of Theorem~\ref{mainthm} also yelds the following result.

\begin{theorem}\label{three}
Let $\Gamma\subset \mathbb{R}^3$ be a $C^{2,\alpha}$ connected simple closed curve such that $T(\Gamma)<4\pi$. There exists $\rho>0$ and $C(\Gamma, \rho)$ such that if $\Gamma'$ is a simple closed curve and $\|\Gamma-\Gamma'\|_{C^{2,\alpha}}\leq\rho$, then 
$$\sup_{\Sigma \in M(\gamma)} |K_\S|\leq C(\Gamma,\rho).$$
\end{theorem}

Let us denote by $E(\Gamma)$ the radius of the largest embedded tubular neighborhood around $\Gamma$. Using Theorem~\ref{three} we can now prove a finer quantative version of Theorem~\ref{mainthm}. For simplicity, we state the theorem assuming that the length of $\Gamma$ is less than  one.

\begin{theorem}\label{mainthm3}
Given $\varepsilon>0, \Delta>0$ and $\theta<4\pi$, there exists a constant $C(\epsilon,\Delta, \theta)$ such that the following holds. If $\Gamma\subset\mathbb{R}^3$ is a $C^{2,\alpha}$ connected simple closed curve whose lenght is less than one, $E(\Gamma)\geq\varepsilon$, $\|\Gamma\|_{C^{2,\alpha}}\leq\Delta$, and $T(\Gamma)\leq \theta$, then
$$\sup_{\S\in M(\Gamma)}|K_\S|\leq C(\epsilon,\Delta, \theta).$$
\end{theorem}

\begin{proof}
Suppose that there exist sequences $\Gamma_n$ and $\S_n\in M(\Gamma_n)$ for which the curvature goes to infinity. The conditions on $\Gamma_n$ guarantee that there exists a connected simple closed curve $\Gamma$ whose total curvature is less than or equal to $\theta$, and a subsequence $\Gamma_{n_k}$ such that $\|\Gamma-\Gamma_{n_k}\|_{C^{2,\alpha}}$ is going to zero. We can then apply Theorem~\ref{three} to reach a contradiction.
\end{proof}
 
\section{Applications}

The results below are consequences of the theorems discussed in the previous sections.

\begin{corr}
If $\Gamma\subset \mathbb{R}^3$ is a $C^{2,\alpha}$, connected, simple closed curve such that $T(\Gamma)< 4\pi$, then  $M(\Gamma)$ is compact.
\end{corr}

\begin{corr}
Let $\Gamma\subset \mathbb{R}^3$ be a $C^{2,\alpha}$, connected, simple closed curve such that $T(\Gamma)< 4\pi$. There exists a constant $N(\Gamma)$ such that the genus of any $\S\in M(\Gamma)$ is less than or equal to $N(\Gamma)$.
\end{corr}

\begin{proof}
The Gauss-Bonnet Theorem states that $$\int_{\Gamma} \vec{k}\cdot\vec{n}ds+\int_{\Sigma}K_\Sigma=2\pi\chi(\S),$$ where $\vec{k}$ is the curvature vector of the curve $\Gamma$, $\vec{n}$ is the exterior normal of $\S$, and $\chi(\S)$ is the Euler characteristic of $\S$. The first integral is bounded in absolute value by the total curvature of $\Gamma$, while the second integral is bounded by the area of $\S$ times the bound on the curvature given by Theorem~\ref{mainthm}. Hence, $|\chi(\S)|$ is bounded and so is the genus.
\end{proof}

\begin{theorem}\label{maincor}
The set of $C^{2,\alpha}$, connected, simple closed curves with $T(\Gamma)<4\pi$ and which do not bound an embedded and orientable minimal surface of genus greater than $g$, for any given $g$, is open in the $C^{2,\alpha}$ topology.
\end{theorem} 
\begin{proof}
Let $\Gamma$ and $N(\Gamma)$ be as above. It suffices to show that for any $C^{2,\alpha}$, connected, simple closed curve $\Gamma$  there exists an $\epsilon>0$ such that if $\Gamma'$ is a connected, simple closed curve such that $\|\Gamma-\Gamma'\|_{C^{2,\alpha}}<\epsilon$, then the genus of any $\S\in M(\Gamma')$ is bounded by $N(\Gamma)$. Assuming that the statement is false, for any $n>0$ there exist  $\Gamma_n$ and $\S_n\in M(\Gamma_n)$ such that $\|\Gamma-\Gamma_n\|_{C^{2,\alpha}}<\frac1n$ and  the genus of $\S_n$ is greater than $N(\Gamma)$. However, the second fundamental forms of the $\S_n$'s are uniformly bounded. Therefore, after going to a subsequence, $\S_n$ converges to a surface $\S\in M(\Gamma)$ with genus greater than $N(\Gamma)$, which is a contradiction.
\end{proof}

Using a compactness argument like the one in the proof of Theorem~\ref{mainthm}, one can prove that the bound on the genus depends on the constants described in Theorem~\ref{mainthm}:

\begin{theorem}
Given $\varepsilon>0, \Delta>0$ and $\theta<4\pi$ there exists a constant $N(\epsilon,\Delta, \theta)$ such that the following holds. Let $\Gamma\subset\mathbb{R}^3$ be a $C^{2,\alpha}$ connected simple closed curve whose lenght is less than one and such that $E(\Gamma)\geq\varepsilon$, $\|\Gamma\|_{C^{2,\alpha}}\leq\Delta$, and $T(\Gamma)\leq \theta$, then the genus of any $\S\in M(\Gamma)$ is bounded by $N(\Gamma)$.\end{theorem}

\vspace{1cm}

\center{Giuseppe Tinaglia gtinagli@nd.edu   \\
Mathematics Department, University of Notre Dame, Notre Dame, IN,
46556-4618}
\bibliographystyle{plain}
\bibliography{mybib}
\end{document}